\documentclass{amsart}
\usepackage{epsf}
\def\hepsffile{\leavevmode\epsffile}

\theoremstyle{plain}
\newtheorem{thm}{Theorem}[subsection]

\newtheorem{lem}[thm]{Lemma}

\newtheorem{prop}[thm]{Proposition}

\theoremstyle{definition}
\newtheorem{defin}[thm]{Definition}

\newtheorem{emf}[thm]{}
\newtheorem{rem}[thm]{Remark}

\def\Im{\protect\operatorname{Im}}

\def\Int{\protect\operatorname{Int}}

\def\mod{\protect\operatorname{mod}}

\def\pr{\protect\operatorname{pr}}

\def\C{{\mathbb C}}
\def\Z{{\mathbb Z}}
\def\Q{{\mathbb Q}}
\def\R{{\mathbb R}}
\def\N{{\mathbb N}}

\def\1{\hbox{\rm\rlap {1}\hskip.03in{\rom I}}}
\def\Bbbone{{\rm1\mathchoice{\kern-0.25em}{\kern-0.25em}
	{\kern-0.2em}{\kern-0.2em}I}}

\def\p{\partial}

\def\wt{\widetilde}
\def\pp{\medskip{\parindent 0pt \it Proof.\ }}

\begin{document}
\hyphenation{Ca-m-po}
\title[Relative Framing of Transverse knots]
{Relative framing of transverse knots}
\author[V.~Chernov (Tchernov)]{Vladimir Chernov (Tchernov)}
\address{Department of Mathematics, 6188 Bradley Hall, Dartmouth College, Hanover NH 03755, USA}
\email{Vladimir.Chernov@dartmouth.edu}
\begin{abstract}
It is well-known that a knot in a contact manifold $(M,C)$ transverse to a trivialized contact structure possesses the natural framing given by the first of the trivialization vectors along the knot. If the Euler class $e_C\in H^2(M)$ of $C$ is nonzero, then $C$ is nontrvivializable and the natural 
framing of transverse knots does not exist.

We construct a new framing-type invariant of transverse knots called relative framing.
It is defined for all tight $C$, all closed irreducible atoroidal $M$, and many other cases when $e_C\neq 0$ and the classical framing invariant is not defined. We show that the relative framing distinguishes many transverse knots that are isotopic as unframed knots. It also allows one to define the Bennequin invariant of a zero-homologous transverse knot in contact manifolds $(M,C)$ with $e_C\neq 0$.

Our recent result is that the groups of Vassiliev-Goussarov invariants of transverse and of framed knots are canonically isomorphic, when $C$ is trivialized and transverse knots have the natural framing. We show that the same result is true whenever the relative framing is well-defined.
\end{abstract}



\maketitle

\section{Introduction} 
We work in the smooth category. 
Throughout this paper $M$ is a (not necessarily compact) 
oriented $3$-dimensional Riemannian manifold, and $C$ is an oriented
$2$-dimensional subbundle of $TM$.

$C$ is called {\em a contact structure\/} if it can be locally presented as $\ker \alpha$, for some $1$-form $\alpha$ 
with $\alpha \wedge d \alpha\neq 0$. All contact manifolds $(M,C)$ have the natural orientation given by $\alpha \wedge d\alpha$. The contact Darboux Theorem says that every contact $3$-manifold is locally contactomorphic to $\bigl (\R^3=(x,y,z), \ker (ydx-dz)\bigr )$. The charts where 
$(M,C)$ is contactomorphic to $\bigl (\R^3=(x,y,z), \ker (ydx-dz)\bigr )$ are called {\em Darboux charts.\/}
If $C$ can be globally presented as $\ker \alpha$ with  $\alpha \wedge d \alpha\neq 0$, then it is called {\em a cooriented (transversally oriented) contact structure.\/} 
In this paper all the contact structures are assumed to cooriented.

A {\em curve\/} in $M$ is an immersion of $S^1$ into
$M$. A {\em framed curve\/} in $M$ is a curve equipped with a continuous unit
normal vector field.
A {\em transverse curve\/} in $(M,C)$
is a curve that is nowhere tangent to the planes of $C$. 
In the case where $C$ is a contact structure the curves that
are everywhere tangent to $C$ are called {\em Legendrian.\/} 
{\em A homotopy\/} of ordinary, framed, transverse, or Legendrian curves is a path 
in the corresponding space of curves.

A {\em knot (resp. framed knot, resp. transverse knot)\/} 
in $(M,C)$ 
is an embedding (resp. framed embedding, resp. transverse embedding) of $S^1$. In the case where $C$ is contact we 
define {\em Legendrian knots\/} as Legendrian curves that are embeddings. An {\em isotopy\/}
of ordinary, framed, transverse, or Legendrian knots is a path in the corresponding 
space of curves that consists of smooth embeddings.
Ordinary, framed, transverse, and Legendrian knots are studied up to the corresponding
isotopy equivalence relation.

It is a classical fact that 
if the $2$-plane bundle $C$ is trivialized, then every transverse
curve $\gamma$ has a natural framing given by the unit normals corresponding to
the projections  of the first vectors of the trivialization of $C$ along $\gamma$
to the $2$-planes orthogonal to the velocity vectors of $\gamma$.

Clearly, if two transverse knots are not homotopic, then they are also not isotopic.
Thus the natural question of the
transverse knot theory in contact manifolds with trivialized $C$ 
is to find and study 
examples of (transverse homotopic) transverse knots that are not isotopic as transverse knots even though 
they realize isotopic framed knots. First such examples of transverse knots in the standard contact $\R^3$ were constructed recently in the ground breaking works of J.~Birman and W.~Menasco~\cite{BirmanMenasco2} and of J.~Etnyre and K.~Honda~\cite{EtnyreHonda}.

What should transverse knot theory study when $C$ is not
trivializable? On the first glance, since there
is no way to assign a framing to a transverse knot, it should study examples of 
(transverse homotopic) transverse knots that are not isotopic as transverse knots even though they
are isotopic as unframed knots.

We show 
that, in a sense, the natural framing of transverse curves 
is still well-defined for a vast collection of $(M,C)$ with  non-trivializable $C$, see
Theorem~\ref{main}.
This means that after we choose a framing 
for a transverse knot $K$, we get the natural mapping from transverse to
framed isotopy classes of knots for all transverse knots that
are homotopic to $K$ as transverse curves. We call such a mapping {\em a relative framing.\/}

Thus for the cases where relative framing is defined 
the natural question of the transverse knot theory 
should be once again to distinguish
transverse knots that realize isotopic framed knots.
In particular, when $C$ is contact we construct a vast collection of transverse knots that are different as transverse knots but are not distinguishable by the previously known invariants. The relative framing, when it is defined, provides an alternative definition of the Bennequin invariant of zero-homologous transverse knots in 
contact $(M,C)$ with non-trivializable $C$, see~\ref{Bennequin}. The previously known definitions of the Bennequin invariant in the case of nontrivializable $C$, see for example~\cite{EtnyreLegendrianTransversal},  were dependent on the choice of the relative homology class in the knot complement. The only ambiguity in our construction of the Bennequin invariant is the choice of an additive constant for  each homotopy class of transverse curves.

A recent result of the author~\cite{ChernovLegendrian} is that the groups of Vassiliev-Goussarov invariants of framed knots and of knots transverse to a contact structure $C$ are canonically isomorphic when $C$ is trivialized, and hence transverse knots have the well-defined natural framing. This appears to be also true whenever relative framing of transverse knots is well-defined.  

\section{Relative framing}
Everywhere below we denote by $\mathcal C$ a connected component of the space of unframed curves in $M$. (As it follows from the Smale-Hirsch $h$-principle, connected components of the space of curves in $M$ are naturally identified with the conjugacy classes of $\pi_1(M)$, see~\ref{componentstransverse}.) 

\subsection{Covering $\pr:\mathcal F\to \mathcal C$}\label{subsectioncovering} (In the presentation of the covering $\pr:\mathcal F\to \mathcal C$, that was introduced in the first version of this text we follow our work~\cite{ChernovFramed}.)

Let $\mathcal F'$ be the subspace of the space of framed curves in
$M$ that consists of all framed curves realizing curves from $\mathcal
C$ if we forget the framing. (In fact it is possible to show that this subspace consists of two connected components of the space of framed curves in $M$.)
Put $\pr':\mathcal F'\rightarrow
\mathcal C$ to be the forgetting of the framing mapping.

Let $p:\mathcal F'\to \mathcal F$ be the quotient by the 
following equivalence relation: $f'_1\sim f'_2$ if 
there exists a path $q:[0,1]\rightarrow \mathcal F'$ connecting $f'_1$ and
$f'_2$ such that $\Im(\pr'(q))=\pr'(f'_1)=\pr'(f'_2)$. (This means 
that we
identify two framed curves if the nonzero sections of the normal bundle to the curve induced by the framings are homotopic as nonzero sections.)
Put $\pr:\mathcal F\rightarrow \mathcal C$ to be the mapping such that
$\pr\circ p=\pr'$.

\begin{lem}[see~\cite{ChernovFramed}]\label{covering}
$\pr:\mathcal F\rightarrow \mathcal C$ is a regular covering with a structure group $\Z$. 

The  mapping $\delta:\pi_1(\mathcal C, c)\rightarrow \Z$, that maps the class $[\alpha]\in\pi_1(\mathcal C, c)$ of a loop $\alpha:[0,1]\to \mathcal C$ to
the element $\delta([\alpha])$ of the structure group $\Z$ of $\pr$ such that $\delta([\alpha]) \cdot \widetilde \alpha(0)=\widetilde \alpha(1)$, is a homomorphism. (Here $\cdot$ denotes the action of the structure group $\Z$ and $\widetilde \alpha:[0,1]\to \mathcal F$ is a lift of $\alpha$.) Since $\Z$ is abelian, $\delta$ can be also regarded as a homomorphism $\delta:H_1(\mathcal C)\to \Z$.
\end{lem}

The proof of the Lemma is straightforward.

\begin{defin}[of isotopic knots from $\mathcal F$, see~\cite{ChernovFramed}]\label{isotopy}
Let $K_0, K_1\in\mathcal F$ be such that
$\pr(K_0)$ and
$\pr(K_1)$ are knots (embedded curves). Then $K_0$ and $K_1$ are said to be
{\em isotopic\/} if there exists a path $q:[0,1]\rightarrow \mathcal F$ such that 
$q(0)=K_0, q(1)=K_1$, and $\pr\circ q$ is an isotopy of unframed knots.
Lemma~\ref{covering} implies that framed knots $K_{f,0},  K_{f,1}\in\mathcal F'$
are framed isotopic if and only if $p(K_{f,0})$ and
$p(K_{f,1})$ are isotopic in $\mathcal F$.
\end{defin}

The following Theorem~\ref{framedtheorem} that was proved by us in~\cite{ChernovFramed} was first stated 
(for compact manifolds and in a different formulation) by Hoste and Przytycki~\cite{HostePrzytycki}. They referred to the work~\cite{McCullough} of McCullough on mapping 
class groups of 3-manifolds for the idea of the proof of 
this fact. However to the best of our knowledge the proof of this fundamental fact was not given in literature. The proof we provide in~\cite{ChernovFramed}, see Theorem~2.0.5, is 
based on the ideas and methods different from the ones 
Hoste and Przytycki had in mind. (In the case of $[K]=0\in H_1(M)$, the statement of Theorem~\ref{framedtheorem} is obvious because the self-linking invariant of a zero-homologous framed knot is well-defined.)

\begin{thm}\label{framedtheorem}
Let $M$ be a manifold that does not contain an embedded non-separating $S^2$ (or equivalently $M$ is not realizable as a connected sum $M=(S^1\times S^2)\#M'$). 
Then $i\cdot K$ is not isotopic in $\mathcal F$ to $K$ for any nonsingular knot $K$ and any $i\neq 0$. (Here $\cdot $ denotes the action of the structure group $\Z$ of the regular covering $\pr :\mathcal F\to \mathcal C$.)
\end{thm}

\subsection{Covering $\pr_{\mathcal T}:\mathcal F_{\mathcal T} \to \mathcal T$ and relative framing}
Let $(M,C)$ be an oriented $3$-manifold equipped with an oriented 
$2$-plane bundle $C$ and let $\mathcal T\subset \mathcal C$ be a connected component of the space of transverse to $C$ curves.  Let $\mathcal F_{\mathcal T}$ be the space that consists
of pairs $(\tau, \phi)$, where $(\tau:S^1\to M)\in \mathcal T$ and $\phi$ is a homotopy class of a nowhere zero section of $\tau^*(C)$.

\begin{lem}\label{trivialization}
The natural mapping $\pr_{\mathcal T} : \mathcal F_{\mathcal T} \to \mathcal T$ is a regular covering with a structure group $\Z$. 
Let $e_C\in H^2(M)$ be the Euler class of $C$, let $[S^1\times S^1]\in H_2(S^1\times S^1)$ be the fundamental class of $S^1\times S^1$, and let 
$m_{\mathcal T}=GCD\Big \{ f^*(e_C)([S^1\times S^1])\in \Z \Big | f:S^1\times S^1\to M \text{ such that } f\big |_{S^1\times s_0}\in \mathcal T, \text{ for all } s_0\in S^1\Big \}.$
Then the stabilizer of the set of connected components of $\mathcal F_{\mathcal T}$ is a subgroup $m_{\mathcal T}\Z \subset \Z$. 
\end{lem}

The proof of this Lemma is straightforward. To get that $m_{\mathcal T}\Z\subset \Z$ is the stabilizer one uses the obstruction theory interpretations of the Euler class and of the relative Euler class. 
(Recall that the relative Euler class is the obstruction for the existence 
of a nonzero extension of a nonzero section of a bundle given over manifold's boundary $\p N$ to the whole manifold $N$.)  

\begin{defin}[of relative framing]
For transverse curves in $(M,C)$ the normal bundle to a curve $(\tau: S^1\to M)\in\mathcal T$ can be canonically identified with the pull back $\tau^*(C)$ of $C$. Thus 
$\pr_{\mathcal T}:\mathcal F_{\mathcal T}\to \mathcal T$ is the restriction of $\pr : \mathcal F\to \mathcal C$ to $\mathcal T\subset \mathcal C$. 

By Lemma~\ref{trivialization} the quotient covering 
$\wt { \pr_{\mathcal T}}:\mathcal F_{\mathcal T}/(m_{\mathcal T}\Z) \to \mathcal T$ has a section 
$F_{m_{\mathcal T}}:\mathcal T\to \mathcal F_{\mathcal T}/(m_{\mathcal T}\Z)$ that is uniquely defined modulo the choice of its value on a point of $\mathcal T$.
{\em We call such a section a $\mod m_{\mathcal T}$ relative framing.\/} In particular, if $m_{\mathcal T}=0$, then  
$\pr_{\mathcal T}:\mathcal F_{\mathcal T} \to \mathcal T$ has a section $F:\mathcal T\to \mathcal {F_T}$ and we call this section  a {\em  relative framing.\/}  We call 
$c\in \mathcal F/ (m_{\mathcal T}\Z)$ a {\em knot (in $\mathcal F/ (m_{\mathcal T}\Z)$)\/} if $\wt {\pr_{\mathcal T}}(c)\in\mathcal C$ is an (embedded) knot. We define {\em an isotopy in $\mathcal F/ (m_{\mathcal T}\Z)$\/} to be a path that projects to an isotopy in $\mathcal C$. 

The relative framing is a very powerful transverse knot invariant, since if for $K_1, K_2\in \mathcal T$ and a relative framing $F_{m_{\mathcal T}}$, the knots $F_{m_{\mathcal T}}(K_1)$ and $F_{m_{\mathcal T}}(K_2)$ are not isotopic in $\mathcal F/ (m_{\mathcal T}\Z)$, then $K_1$ and $K_2$ are clearly not isotopic as transverse knots from $\mathcal T$. 
\end{defin}

\subsection{Examples of relative framings}
The classical framing invariant of transverse to $C$ knots is $(M,C)$ is defined for trivialized $C$. Thus it is impossible to define 
when the Euler class $e_C$ of $C$ is nonzero, and hence $C$ is nontrivializable.
Below we give many examples of $(M,C)$ with $e_C\neq 0$ for which the relative framing exists for all the connected components $\mathcal T$ of the space of transverse to $C$ curves in $(M,C)$.

Recall a few basic definitions. A contact structure is said to be
{\em overtwisted\/} if there exists a $2$-disk $D$ embedded
into $M$ such that the boundary $\p D$ is tangent to $C$ while the disk $D$
is transverse to $C$ along $\p D$. Not overtwisted contact structures are
called {\em tight\/}. Overtwisted contact structures are easy to construct, see Lutz~\cite{Lutz}. A result of Eliashberg~\cite{Eliashberg} says that
every oriented $2$-plane distribution in $TM$ is homotopic to an overtwisted
contact structure. Tight contact structures on the contrary are quite hard
to construct and their classification is one of the main current goals of
the $3$-dimensional contact topology.

A manifold $M$ is said to be {\em irreducible\/} if every embedded
$2$-sphere in $M$ bounds a ball. A closed orientable $3$-manifold is said to be {\em atoroidal \/} 
if for every
mapping $\mu:T^2=S^1\times S^1\rightarrow M$ 
the kernel of $\mu_*:\pi_1(T^2)\rightarrow \pi_1(M)$ is nontrivial.

\begin{thm}\label{main}
Let $(M,C)$ be such that at least one of the following three conditions hold: 
\begin{description}
\item[1] $e_C\in H^2(M)$ is an element of finite order ($C$ is not necessarily contact);
\item[2] $M$ is closed irreducible and atoroidal, in particular if $M$ is closed and admits a Riemannian metric of negative sectional curvature ($C$ is not necessarily contact);
\item[3] $C$ is a tight contact structure;
\end{description} 
then for all the connected components $\mathcal C$ and all $\mathcal T\subset \mathcal C$ there exists a relative framing $F:\mathcal T\to \mathcal F_{\mathcal T}$.

\end{thm}

\begin{rem} 
Using the results of Eliashberg~\cite{Eliashberg} and Lutz~\cite{Lutz}  it is rather easy to show that for every $\alpha\in H^2(M)$ there exists an overtwisted 
cooriented contact structure $C$ on $M$ with $e_C=2\alpha$, see for example~\cite{ChernovLegendrian}. This allows one to construct many examples of contact manifolds $(M,C)$ that satisfy  all the conditions of Theorem~\ref{examples}.
\end{rem}

\pp The proof of Theorem~\ref{main} is based on the observation that it suffices to show that  
$f^*(e_C)([S^1\times S^1])=0$, for all mappings $f:S^1\times S^1\to M$.  Then it is automatically true for those $f$ that satisfy an extra  condition $f\big |_{S^1\times s_0}\in \mathcal T$, for all $s_0\in S^1$. Hence by Lemma~\ref{trivialization} $m_{\mathcal T}=0$, for all $\mathcal C$ and $\mathcal T\subset \mathcal C$.

If $e_C$ is an element of finite order, then $e_C(\alpha)=0$ for all 
$\alpha\in H_2(M, \Z)$. Thus $e_C(f_*([S^1\times S^1]))=f^*(e_C)([S^1\times S^1])=0$, for all $f:S^1\times S^1\to M$, {\em and this proves the Theorem in the case of  condition 1.\/}

The Theorem of Gabai (see~\cite{Gabai} Corollary $6.18$) implies that every
$\alpha\in H_2(M)$ that is realizable by a  mapping of $S^1\times S^1$ can be realized by a collection of embedded spheres and tori. Eliashberg~\cite{Eliashbergtight} proved that $e_C$ vanishes on all embedded tori and spheres, provided that $C$ is tight. Combining
the results of Gabai and Eliashberg we get that $f^*(e_C)([S^1\times S^1])=0$, for all $f:S^1\times S^1\to M$ (and not only for embeddings).
{\em This proves the Theorem in the case of  condition 3.\/}

Since $e_C$ takes zero value on any surface realizing a homology class that is
in the torsion of $H_2(M)$, we get that {\em to prove the Theorem in the case of  condition 2\/} it suffices to show
that if $M$ is closed irreducible and atoroidal, 
then $f_*([S^1\times S^1])$ is in the torsion of $H_2(M)$, for all $f:S^1\times S^1\to M$. 
(In fact,  since $H_2(M)=H^1(M)$ by the Poincare duality and $H_1(M)$ is torsion free, this implies that $f_*([S^1\times S^1])=0\in H_2(M)$.) 
Thus to prove the Theorem it suffices to show that for every $f:S^1\times S^1\to M$ there exists a finite covering $p:S^1\times S^1\to S^1\times S^1$ such that $f_*\circ p_*([S^1\times S^1])=0\in H_2(M)$.

The Sphere Theorem, see for example~\cite{Hempel}, 
says that $\pi_2(M)=0$ for irreducible $M$, and the elementary obstruction theory implies that homotopy classes of maps $f:S^1\times S^1\to M$ are classified by the homomorphisms $f_*:\pi_1(S^1\times S^1)\to \pi_1(M)$.

Let $m,l\in\pi_1(S^1\times S^1, 1\times 1)$ be
the classes of the meridian $S^1\times 1$ and the longitude
$1\times S^1$
of $S^1\times S^1=T^2$. (Here we regard $S^1$ as $\{z\in\C\big | |z|=1\}$.) Since $M$ is atoroidal, we get that
there exist $i,j\in\Z$, with at least one of them nonzero, such that 
$f_*(m)^if_*(l)^j=1\in\pi_1(M)$. 

{\em Assume that both $i$ and $j$ are nonzero.\/} Let $p:S^1\times
S^1\rightarrow S^1\times S^1$ be the covering such that $p_*(m)=m^i$ and
$p_*(l)=l^{-j}$. Clearly $f_*\circ p_*(m)=f_*\circ p_*(l)$, and then $f\circ p$ is homotopic to a mapping that passes through some $\gamma:S^1\to M$. Thus $f_*\circ p_*([S^1\times S^1])=0\in H_2(M)$ and $f_*([S^1\times S^1])\in H_2(M)$ is an element of finite order.

{\em Consider the case where $i=0$ and $j\neq0$.\/} (The case of $j=0$ and
$i\neq 0$ is treated in the same way.) Let $p:S^1\times
S^1\rightarrow S^1\times S^1$ be the covering such that $p_*(m)=m$ and
$p_*(l)=l^j$. Clearly $f_*\circ p_* (l)=1\in \pi_1(M)$, and then $f\circ p$ is homotopic to a mapping that passes through a projection of $S^1\times S^1$ to a meridian. Thus $f_*\circ p_*([S^1\times S^1])=0\in H_2(M)$ and $f_*([S^1\times S^1])\in H_2(M)$ is an element of finite order.  

The fact that every closed manifold $M$ with a Riemannian metric of negative sectional curvature is irreducible follows from the Hadamard and the Sphere Theorems, see for example~\cite{DoCarmo} and~\cite{Hempel}. The fact that such $M$ is atoroidal is an immediate consequence of the Preissman Theorem, see for example~\cite{DoCarmo}. 
\qed

\begin{emf}\label{Bennequin}{\em Bennequin invariants for transverse knots in tight contact $(M,C)$ with non-trivializable $C$.\/}
Clearly the relative framing can be used to define the self-linking number of a zero-homologous transverse knot $K$ which is the famous Bennequin invariant $\beta(K)$. (The homology definition of the self-linking invariant does not work unless $[K]=0\in H_1(M)$, but it is possible to define various affine generalizations of the self-linking invariant for nonzero homologous framed knots, see our work~\cite{ChernovFramed}.)
The Bennequin invariant $\beta=\beta_{F}$ defined this way depends
on the choice of the relative framing $F:\mathcal T\to \mathcal F$ used to define $\beta$, but $\beta_{F_1}-\beta_{F_2}$ is a constant function on transverse knot isotopy classes of transverse knots from $\mathcal T$, for any two relative framings $F_1, F_2:\mathcal T \to \mathcal F$.

The previously known way, see for example~\cite{EtnyreLegendrianTransversal}, to define the Bennequin invariant of a zero homologous transverse knot $K$ in $(M,C)$ with $e_C\neq 0$ is described below. Take a Seifert surface $\Sigma_K$ of $K$.  Define the self-linking number of $K$ using the framing of $K$ given by vectors in contact planes along $C$ such that the section of $C$ over $\p \Sigma_K=K$ given by these vectors extends to a nonzero section of $C$ over $\Sigma_K$. Since $e_C\neq 0$, it is clear that for a transverse knot $K$ this definition of the Bennequin invariant of a knot $K$ does depend on the choice of the relative homology class $[\Sigma_K]\in H_2(M\setminus \Int T_K , \p T_K)$. 
(Here $T_K$ is the tubular neighborhood of $K$ in $M$.)
The other drawback of this definition is that it is not clear how to identify relative homology classes for two different transverse knots, and thus it is not clear how to appropriately compare the values of Bennequin invariants for two different transverse knots.

Our definition of the Bennequin invariant has much less ambiguity and allows one to compare the values of Bennequin invariants for two transverse knots inside $\mathcal T$.

\end{emf}

\subsection{Some useful facts about transverse curves and relative framings}
\begin{emf}\label{h-principletransverse}
{\em $h$-principle for transverse curves.\/} For a topological space $X$ we denote by $\Omega_X$ the space of free loops in $X$.

Let $(M,C)$ 
be a contact manifold and let $TM\setminus C$ be the topological space that
is $TM$ with the contact subbundle cut out. 
The $h$-principle for transverse curves proved by M.~Gromov~\cite{Gromov}
p.84, see also~\cite{EliashbergMishachev} Section 14.2, says 
that the space of transverse curves in $(M,C)$ is weak
homotopy equivalent to $\Omega_{TM\setminus C}$. The equivalence is given by mapping a point $t$ of $S^1$ that
parameterizes a transverse curve $\tau$ to the point of $TM\setminus C$ that corresponds to the velocity vector of $\tau$ at $\tau(t)$.

If the contact structure is cooriented, then $TM\setminus C$ is homotopy
equivalent to $M\times S^0$, and thus the space of transverse curves in
such $(M,C)$ is weak homotopy equivalent to $\Omega_
{M\times S^0}$. Identify $S^0$ with $\{-1,1\}$. Then the
equivalence is given by mapping a point $t$ of $S^1$ that 
parameterizes a transverse curve $\tau$ to $\tau(t)\times \{1\}$, provided that 
the velocity vectors of $\tau$ point in the direction of coorienting
half-spaces (the ones where the $1$-form $\alpha$ from the definition of a cooriented contact structure is positive); and by mapping it to $\tau(t)\times \{-1\}$, otherwise.
\end{emf}

\begin{emf}\label{h-principleimmersed} {\em The Hirsch-Smale $h$-principle for immersed curves,\/} see for example~\cite{Gromov}, says that the space of immersed curves in $M$ is weak homotopy equivalent to $\Omega_{STM}$, where $STM$ is the total space of the unit tangent bundle $STM\to M$. From the homotopy sequence of the locally trivial 
$S^2$-fibration $STM\to M$ we get that the set of connected components of the space of curves in $M$ is identified with the set of connected components of $\Omega_M$.
\end{emf}

Combining the $h$-principles for curves and transverse curves we get the following proposition.

\begin{prop}\label{componentstransverse}
The set of connected components of the space of curves in $M$ is naturally identified with the conjugacy classes of $\pi_1(M)$. Let $C$ be a cooriented contact structure on $M$, then every connected component
$\mathcal C$ of the space of curves in $M$ contains precisely two connected components of the space of transverse curves in $(M,C)$. The two connected components are distinguished by whether the orientation of an (immersed) transverse curve points in the direction of the coorienting $C$ half-space of $TM\setminus C$ or not.
\end{prop}

\begin{emf}\label{description} 
{\em Description of transverse curves.\/}
Identify a  point $(x, y, z)\in \R^3$ with the
point $(x,z)\in \R^2$ furnished with the fixed
direction of an unoriented   
straight line through $(x,z)$ with the slope $y$. Then the curve in $\R^3$
can be presented as a one parameter family of points with non-vertical directions in $\R^2$. 

Let $\alpha=ydx -z$ be the standard contact form on $\R^3$. Clearly the curve in $\R^3$ is transverse to the contact structure $\ker \alpha$ if and only if its projection to $(x,z)$-plane is never tangent to the direction field along itself. For example, following~\cite{FuchsandTabachnikov} we observe that the knot in Figure~\ref{directionknot.fig} is not transverse. (Since we project to the $(x,z)$-plane the standard orientation conventions imply that at the crossings 
the strands with a larger value of $y$-coordinate are located below.)

\begin{figure}[htbp]
 \begin{center}
  \epsfxsize 8cm
  \hepsffile{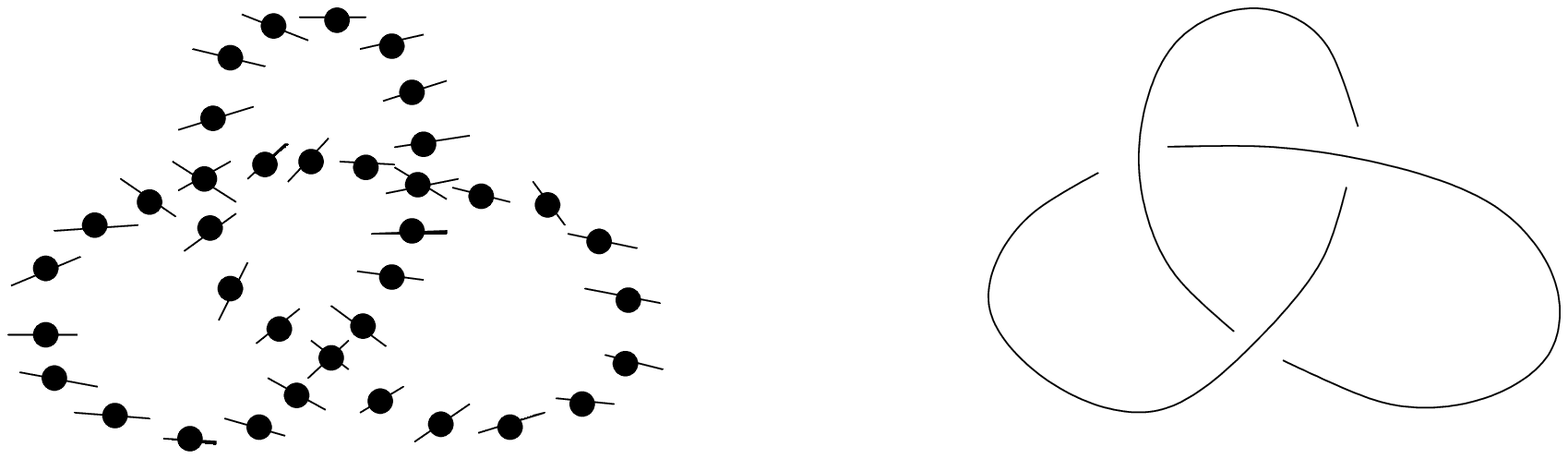}
 \end{center}
\caption{}
\label{directionknot.fig}
\end{figure}

Following Fuchs and Tabachnikov~\cite{FuchsandTabachnikov} we observe that the isotopy class of transverse knots in $(\R^3, \ker \alpha)$, whose velocity vectors point in the direction of $\alpha>0$,
is determined by a standard regular knot diagram with over and under passes, with two extra conditions (determined by the fact that the curve should be transverse to the non-vertical direction field):
\begin{description}
\item[1] the knot diagram should be free of vertical tangents directed upward;
\item[2] if at the crossing the upward vertical direction lies inside of the angle formed by the two oriented tangent lines to the strands, then the strand directed to the upper-right corner is lower than the strand directed to the left-upper corner. 
\end{description}

The two forbidden fragments are shown in Figure~\ref{forbidden.fig}. For example, see~\cite{FuchsandTabachnikov}, the diagram~\ref{trefoil.fig} is a valid diagram of a transverse trefoil knot.

Transverse knots, whose velocity vectors point in the direction of $\alpha<0$, can be described in the similar way with all the orientations of the strands in the forbidden fragments changed to the opposite.

\begin{figure}[htbp]
 \begin{center}
  \epsfxsize 6cm
  \hepsffile{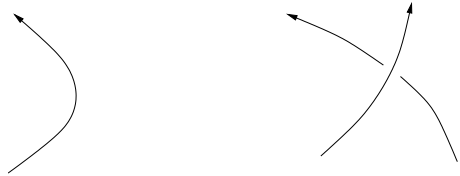}
 \end{center}
\caption{Forbidden fragments of a transverse knot diagram.}
\label{forbidden.fig}
\end{figure}

\begin{figure}[htbp]
 \begin{center}
  \epsfxsize 6cm
  \hepsffile{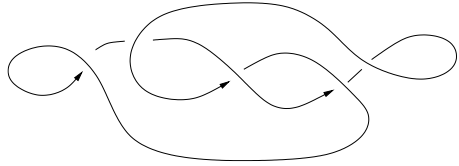}
 \end{center}
\caption{Transverse trefoil knot.}
\label{trefoil.fig}
\end{figure}

\end{emf}

\subsection{Relative framing as an invariant of knots transverse to a contact structure}

The following Theorem provides many examples where relative framing distinguishes transverse knots in contact manifolds $(M,C)$ with $e_C\neq 0$ that are isotopic as unframed knots and are homotopic as transverse curves. (Since $e_C\neq 0$, $C$ is nontrivializable and the classical framing of transverse knots is not defined.) 

\begin{figure}[htbp]
 \begin{center}
  \epsfxsize 8cm
  \hepsffile{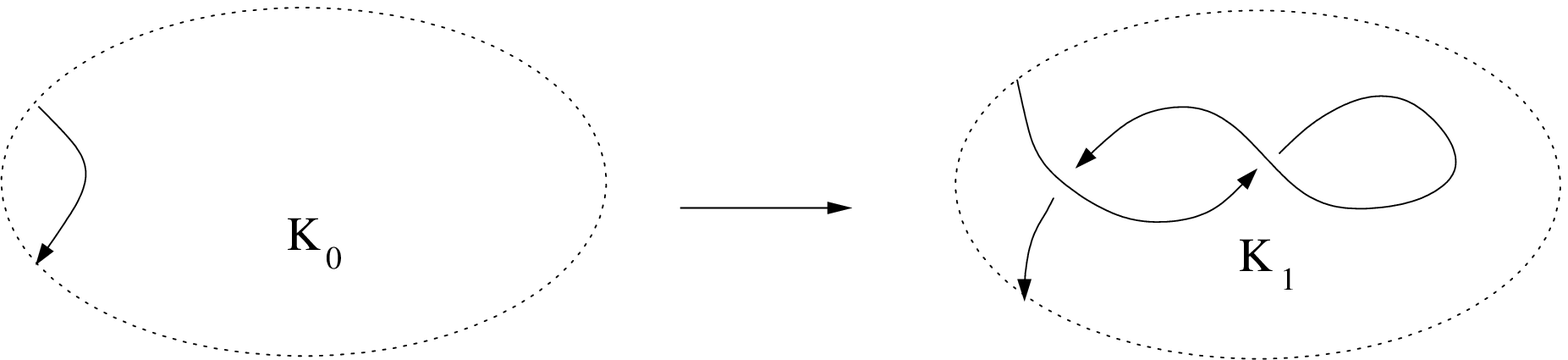}
 \end{center}
\caption{}
\label{stabilization.fig}
\end{figure}

\begin{thm}\label{examples}
Let $(M,C)$ be a contact manifold with a non-trivializable 
cooriented contact structure such that $(M,C)$ satisfies any one of the three conditions of Theorem~\ref{main}, and hence all the connected components of the space of transverse curves in $(M,C)$ admit a relative framing.
 
Let $K_0$ be a transverse knot in $(M,C)$ and let $K_1$ be a transverse
knot that is the same as $K_0$ everywhere except of a small piece located 
in a Darboux chart where 
it is changed as it is shown in Figure~\ref{stabilization.fig} (in terms of the transverse knot diagrams described in~\ref{description}).

Then 
\begin{description}
\item[1]
$K_0$ and $K_1$ realize isotopic unframed knots and belong to the same
component $\mathcal T$ of the space
of transverse curves.
\item[2] The relative framing invariant distinguishes $K_0$ and $K_1$, provided that $K_0$ is zero-homologous or that $M$ is not realizable as a connected sum $(S^1\times S^2)\# M'$.
\end{description}
\end{thm}

\begin{figure}[htbp]
 \begin{center}
  \epsfxsize 9.5cm
  \hepsffile{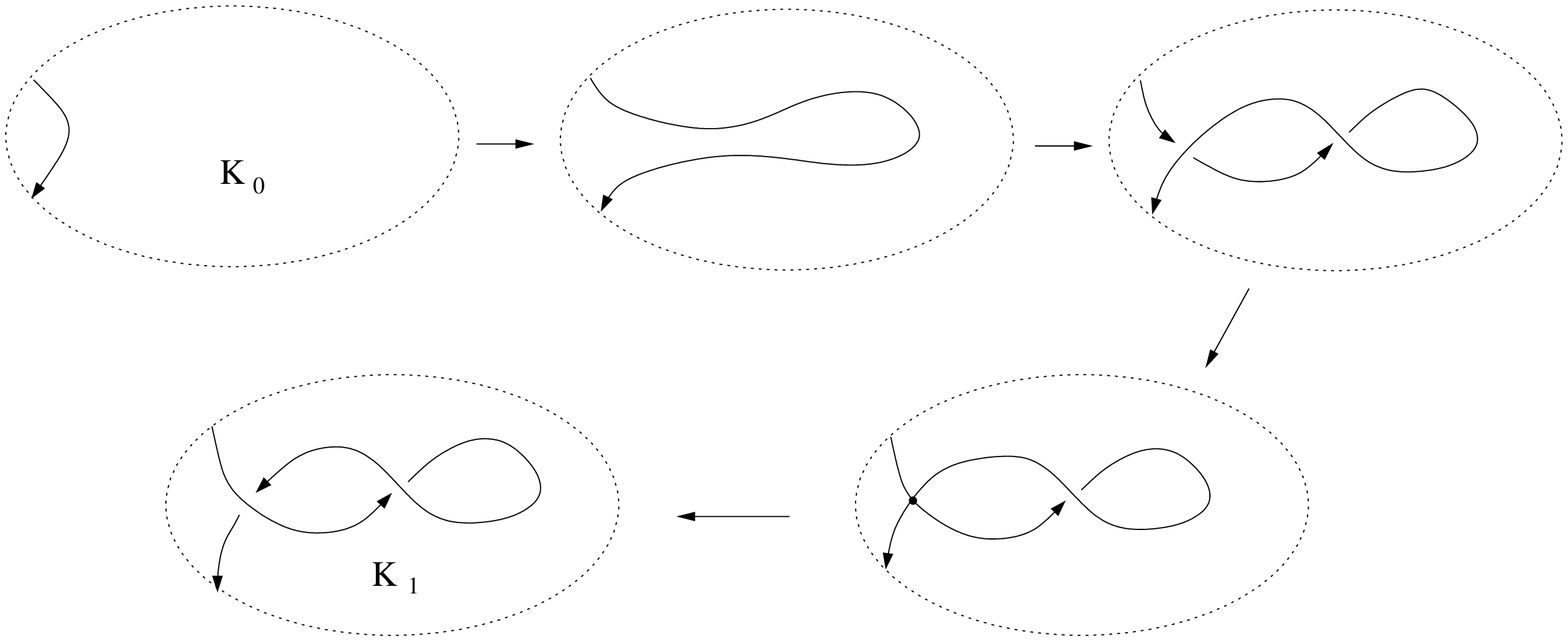}
 \end{center}
\caption{}
\label{kinkcreation.fig}
\end{figure}

\pp Clearly $K_0$ and $K_1$ realize isotopic unframed knots. They belong to the same 
component of the space of transverse curves, since Figure~\ref{kinkcreation.fig} shows  how to deform one of them into the other in the class of transverse curves.

Let $F:\mathcal T\to \mathcal F_{\mathcal T}$ be a relative framing. 
It is easy to see that $F(K_1)$ is isotopic in $\mathcal F$ to $(-2)\cdot F(K_0)$, 
where $\cdot$ denotes the action of the structure group $\Z$ of the covering $\pr:\mathcal F\to \mathcal C$, see Lemma~\ref{covering}. 
If $K_0$ and $K_1$ are isotopic transverse knots, then $F(K_0)$ and $F(K_1)$ are
isotopic in $\mathcal F$. Thus we get that $F(K_0)$ is isotopic to $(-2)\cdot F(K_0)$ in $\mathcal F$. 

For $M\neq (S^1\times S^2)\#M'$ the last statement contradicts to
Theorem~\ref{framedtheorem} and thus $K_0$ and $K_1$ are nonisotopic transverse knots in contact $(M,C)$.

If $K_0$ and $K_1$ are zero homologous knots, then we use the relative framing $F:\mathcal T\to \mathcal F$ to define the Bennequin invariant $\beta_F$ (also known as the self-linking number invariant) of $K_0$ and $K_1$. Since $\beta_F(K_0)\neq \beta_F (K_1)$, we get that $K_0$ and $K_1$ are nonisotopic transverse knots.
\qed

\begin{rem}
Since the contact structure in Theorem~\ref{examples} 
is not trivializable it is impossible to define the classical mapping from transverse to framed knots given by the first vector of trivialization of $C$. Since $K_0$ and $K_1$ are isotopic as unframed knots they are not distinguishable by the classical knot invariants and we see that the relative framing is indeed a very powerful invariant of transverse knots in contact $(M,C)$ with $e_C\neq 0$.
\end{rem}

\section{Vassiliev-Goussarov invariants of knots transverse to a contact structure} Fuchs and Tabachnikov proved~\cite{FuchsandTabachnikov} that for the standard contact $\R^3$ the groups of $\C$-valued Vassiliev invariants of transverse knots and of framed knots are canonically isomorphic. Their proof was based on the existence of the Kontsevich integral~\cite{Kontsevich} universal Vassiliev invariant of framed knots in $\R^3$. The Kontsevich integral was shown to exist in the case of $\C$-, $\R$-, and $\Q$-valued Vassiliev invariants and in the case of $M$ being an oriented total space of an $\R^1$-bundle over an oriented surface $F$ with $\p F\neq \emptyset$ by Andersen, Mattes, Reshetikhin~\cite{AMR} (In the case of $M$ that is the product of an annulus and $\R^1$ the first proof is due to Goryunov~\cite{Goryunov}.) For other manifolds or for Vassiliev invariants with values in other abelian groups the Kontsevich integral is not known to exist.

In~\cite{ChernovLegendrian} we showed that in all the cases where 
the cooriented contact structure is trivialized (and hence all transverse knots in $(M,C)$ have the natural framing) the groups of Vassiliev invariants of transverse and of framed knots with values in any abelian group are canonically isomorphic. (Similar results about Vassiliev invariants of framed, Legendrian, and pseudo-Legendrian knots were obtained by us in~\cite{ChernovLegendrian} and~\cite{ChernovPseudoLegendrian}.)
In this work we show that the groups of Vassiliev invariants of transverse and of framed knots are canonically isomorphic in all the cases when transverse knots have a relative framing.

Below we recall some basic definitions from the theory of Vassiliev invariants. (In our exposition of the definitions we follow~\cite{FuchsandTabachnikov}, \cite{ChernovLegendrian}, and~\cite{ChernovPseudoLegendrian}.)

A {\em singular (framed)\/} knot with $n\geq 0$ double points is a curve (framed curve)
in $M$ whose only singularities are $n$ transverse double points.
An {\em isotopy\/} of a singular (framed) knot 
with $n$ double points is a path in the space of singular (framed) knots with
$n$ double points under which the preimages of the double points on $S^1$
change continuously. In a similar way we define singular transverse knots with $n$ double points in the contact manifolds $(M,C)$ and singular knots with $n$ double points from $\mathcal F$. In the last case an isotopy of a singular knot from $\mathcal F$ is a path in $\mathcal F$ that projects to an isotopy of the corresponding singular knot from $\mathcal C$.

A transverse double point $d$ of a singular knot can be resolved in two 
essentially different ways. A resolution of a double point is
called positive (resp. negative) if the tangent vector to the
first strand, the tangent vector to the second strand, and the vector from
the second strand to the first strand form the positive $3$-frame. (This does 
not depend on the order of the strands).
A singular framed (resp. transverse, resp. from $\mathcal F$) knot $K$ with $(n+1)$ 
transverse double points
admits $2^{n+1}$ possible resolutions of the double points. The sign of the resolution 
is put to be $+$ if the number of negatively resolved double points is even; and
it is put to be $-$, otherwise.

Let $\mathcal A$ be an abelian group, and let $x$ be an $\mathcal A$-valued invariant of framed (resp. transverse, resp from $\mathcal F$) knots. The invariant $x$ is said to be a {\em Vassiliev invariant\/}
if there exists a nonnegative 
integer $n$ such that for any singular knot $K_s$ with $(n+1)$
transverse double points the sum (with appropriate signs) of the values of $x$ on the nonsingular
knots obtained by the $2^{n+1}$ resolutions of the double points is zero. 
{\em (Vassiliev invariants are often also called Vassiliev-Goussarov invariants or finite order invariants.)\/} 
A Vassiliev invariant is said to be of order not greater than $n$ (of order $\leq n$) if $n$
can be chosen as the integer in the definition above. The group of $\mathcal
A$-valued finite order invariants has an increasing filtration by the
subgroups of the invariants of order $\leq n$.

Let $\mathcal T$ be a connected component of the space of transverse curves in $(M,C)$ such that transverse knots from $\mathcal F$ admit 
a relative framing $F:\mathcal T\to \mathcal F_{\mathcal T}\subset \mathcal F$. Let $\mathcal F_0$ be the component of $\mathcal F$ that contains $\Im F$ and let $\mathcal F'_0$ be the corresponding component of the space of framed curves.

One easily proves the following Proposition.

\begin{prop}\label{usefuleasy} 
\begin{description}
\item[1] $q:\mathcal F'_0\to \mathcal F_0$  induces the natural bijection between isotopy classes of singular knots with $n$ double points from $\mathcal F'_0$ and from $\mathcal F_0$, for all $n\geq 0$.
\item[2] The pull back mapping $q^*:V^n_{\mathcal F'_0, \mathcal A}\to V^n_{\mathcal F_0, \mathcal A}$ is the canonical isomorphism of the groups of $\mathcal A$-valued order $\leq n$ Vassiliev invariants of knots from $\mathcal F'_0$ and from $\mathcal F_0$, respectively.
\end{description}
\end{prop}

Clearly the mapping $F:\mathcal T\to \mathcal F$ respects the isotopy classes of transverse knots in a contact $(M,C)$. The following Theorem 
says that the groups of Vassiliev invariants of framed knots from $\mathcal F'_0$ and of transverse knots from $\mathcal T$ are canonically isomorphic.

\begin{thm}\label{Vassiliev}
The natural homomorphism $F^*:V^n_{\mathcal F_0, \mathcal A}\to V^n_{\mathcal T, \mathcal A}$ is a canonical isomorphism (that depends 
on the choice of the relative framing $F:\mathcal T\to \mathcal F_{\mathcal T}$). Hence $F^*\circ q^*:V^n_{\mathcal F'_0, \mathcal A}\to V^n_{\mathcal T, \mathcal A}$ is a canonical isomorphism.
\end{thm}

\begin{emf} {\em Proof of Theorem~\ref{Vassiliev}.\/}
Clearly the pull back of an order $\leq n$ Vassiliev invariant via $F^*$ is a Vassiliev invariant of order $\leq n$. It is also clear that $F^*$ is a homomorphism.

Fuchs and Tabachnikov~\cite{FuchsandTabachnikov} have introduced an operation of addition of a transverse double loop to a transverse knot, shown in Figure~\ref{stabilization.fig}. 
We call an addition of $i$ such double loops to a transverse knot $K$ an {\em $i$-stabilization\/} and denote it by $K^i$. Fuchs and Tabachnikov~\cite{FuchsandTabachnikov} proved that if $K_1, K_2$ 
are transverse knots (in the standard contact $\R^3$) that are isotopic as unframed knots and homotopic as transverse immersions, then there exist $i,j\in \N$ such that $K_1^i$ and $K_2^j$ are isotopic transverse knots. 
They also showed that if $i$ and $j$ can be chosen to be equal, then $\C$-valued Vassiliev invariants of transverse knots do not distinguish $K_1$ and $K_2$. 

As it was later observed by Fuchs and Tabachnikov~\cite{FuchsandTabachnikovprivate} the proof of this statement in fact goes through for Vassiliev invariants with values in any abelian group $\mathcal A$ and for all contact manifolds $(M,C)$ with a cooriented contact structure. The standard contact structure on $\R^3$ is trivialized and the $1$-stabilization decreases the self-linking number of the naturally framed transverse knot by two. Thus if $K_1$ and $K_2$ realize isotopic framed knots in $\R^3$, then $i$ and $j$ automatically have to be equal. If $M\neq \R^3$, then knots are not necessarily zero-homologous and the self-linking invariant of framed knots is not well-defined. 

In~\cite{ChernovLegendrian} we showed that if Vassiliev invariants do not distinguish Legendrian (resp. transverse) knots that are isotopic as framed knots, then the groups of Vassiliev invariants of Legendrian (resp.~transverse) and of framed knots are canonically isomorphic. This statement would have being obvious provided that every framed knot was realizable by a Legendrian (resp. transverse) knot. However the famous Bennequin inequality shows that this is not so even for the standard contact $\R^3$. A straightforward verification shows that the same fact is true for Vassiliev invariants of transverse knots from $\mathcal T$ and of knots from $\mathcal F_0$.

Thus the only thing we have to show is that if transverse knots $K_1, K_2\in \mathcal T$ are such that $F(K_1)$ and $F(K_2)$ are isotopic in $\mathcal F_0$, then $i$ and $j$ such that $K_1^i$ and $K_2^j$ are transverse isotopic can be chosen to be equal.

{\em If $M$ does not contain embedded non-separating spheres\/} (or which is the same $M$ is not realizable as a connected sum $M=(S^1\times S^2)\#M'$), then the proof is very elegant.  Namely since $K_1^i$ and $K_2^j$ are isotopic transverse knots, we get that $F(K_1^i)$ and $F(K_2^j)$ are isotopic in $\mathcal F$. Also it is clear that $F(K_1^i)$ is isotopic (in $\mathcal F$) to $(-2i) \cdot F(K_1)$,
and $F(K_2^j)$ is isotopic (in $\mathcal F$) to $(-2j) \cdot F(K_2)$. (Here $\cdot$ denotes the action of the structure group $\Z$ of the covering $\pr :\mathcal F\to \mathcal C$, see Lemma~\ref{covering}.) Since $\pr:\mathcal F\to \mathcal C$ is a regular covering, we get that $F(K_1)$ is isotopic (in $\mathcal F$) to $2(i-j)F(K_1)$.
Theorem~\ref{framedtheorem} implies that $i=j$ and the proof is finished. 

{\em Below we provide the proof without the assumption that $M\neq (S^1\times S^2)\#M'$.\/}
Let $I:[0,1]\to \mathcal F$ be an isotopy between $F(K_1)$ and $F(K_2)$. One verifies that for $i$ big enough the isotopy $\wt I$ of transverse knots $K_1^i$ and $K_2^j$ constructed by Fuchs and Tabachnikov~\cite{FuchsandTabachnikov} can be chosen so that:
\begin{description}
\item[1] the isotopies $\pr (I)$ and $\wt I$ are $C^0$-close; and moreover
\item[2] for every $t\in [0,1]$ the unframed knot $\pr (I(t))$ is $C^0$-close isotopic to the unframed knot 
$\wt I(t)$ inside of the tubular neighborhood of $\wt I (t)$ which is a thin knotted solid torus $T_t\subset M$.
\end{description}

For two framed knots that are isotopic to each other as unframed knots inside of $T_t$ there is a $\Z$-valued obstruction to be isotopic to each other as framed knots inside of $T_t$. This obstruction is the difference between the self-linking numbers of the two knots in the standard solid torus in $\R^3$ defined via some identification between the knotted solid-torus $T_t\subset M$ and the standard solid torus in $\R^3$. 

One verifies that even though the two self-linking numbers do depend on the choice of an identification between $T_t\subset M$ and the standard solid torus in $\R^3$, the difference of
the two self-linking numbers is well-defined and does not depend on this choice. 
Clearly the $\Z$-valued obstruction for $F(\wt I (t))$ and $I(t)$ to be isotopic in $\mathcal F$ (via an isotopy projecting to a $C^0$-small isotopy in the corresponding solid torus) continuously depends on $t$.

Finally, a straightforward verification shows that $i$-stabilization of a transverse knot $K$ corresponds to the addition of $2i$ negative extra twists to the relative framing of $K$. Thus the obstruction for $F(K^i)$ to be $C^0$-small isotopic to $F(K)$ inside of a thin solid torus in $M$ is equal to $-(2i)$. 
The isotopy of transverse knots $\wt I(t)$ between $K_1^i$ and $K_2^j$ was chosen to be  a $C^0$-small approximation for the 
isotopy $\pr (I(t))$ of unframed knots. Since $I(t)$ is an isotopy of $F(K_1)$ to $F(K_2)$ in $\mathcal F$, 
we have that the obstruction for $F(K_1)$ and $F(K_1^i)$ to be $C^0$-small isotopic in $\mathcal F$  equals to the obstruction for $F(K_2)$ and $F(K_2^j)$ to be $C^0$-small isotopic in $\mathcal F$. Thus $i=j$.  \qed
\end{emf}

{\bf Acknowledgments.}
I am grateful to N.~Mishachev for the discussions about the
$h$-principle. 
The first draft of this  paper was written during my stay at the Institute for Mathematics Z\"urich University. I added results about Vassiliev invariants of transverse knots and substantially revised the text at Dartmouth College, when I was supported by the Dartmouth free term research salary.
I  thank Z\"urich University and Dartmouth College for the excellent working conditions.

\end{document}